\documentclass[12pt,a4paper,dvipsnames]{amsart}
\usepackage[utf8]{inputenc}
\usepackage{tabularx}
\usepackage[table]{xcolor}
\usepackage{amsmath, nicefrac, amsthm, verbatim, amsfonts, amssymb, xcolor, bbm}
\usepackage{graphics, xspace, enumerate, cite}
\usepackage[a4paper,margin=3cm]{geometry}%this is for easy margins
\usepackage{mathrsfs}

\usepackage[bb=boondox]{mathalfa}
\usepackage{tikz}
\usepackage{graphicx}
\usepackage[colorlinks=true,citecolor=green,urlcolor=green,linkcolor=green,bookmarksopen=true,unicode=true,pdffitwindow=true]{hyperref}
\usepackage[english]{babel}
\usepackage[languagenames,fixlanguage]{babelbib}

\usepackage[labelformat=simple]{subcaption}

\usetikzlibrary{decorations.pathreplacing,calligraphy}

\makeatletter
\@namedef{subjclassname@2020}{\textup{2020} Mathematics Subject Classification}
\makeatother

\theoremstyle{plain}
\newtheorem{theorem}{Theorem}[section]

\theoremstyle{definition}
\newtheorem{remark}[theorem]{Remark}

\usepackage{textcomp}

\newcommand{\hull}{\operatorname{hull}}
\newcommand{\arctg}{\operatorname{arctg}}
\newcommand{\arcctg}{\operatorname{arcctg}}
\newcommand{\BM}{\operatorname{BM}}
\newcommand{\BB}{\operatorname{BB}}
\newcommand{\erf}{\operatorname{erf}}

\newcommand{\sfL}{\mathsf{L}}
\newcommand{\sfA}{\mathsf{A}}
\newcommand{\sfT}{\mathsf{T}}
\newcommand{\bfe}{\mathbf{e}}
\newcommand{\calL}{\mathcal{L}}
\newcommand{\calA}{\mathcal{A}}

\newcommand{\cH}{\mathrm{H}}
\newcommand{\fatW}{\mathbf{W}}
\newcommand{\fatB}{\mathbf{B}}

\newcommand{\bbN}{\mathbb{N}}

\newcommand{\bbR}{\mathbb{R}}

\newcommand{\bbP}{\mathbb{P}}
\newcommand{\bbE}{\mathbb{E}}
\newcommand{\bbjedan}{\mathbbm{1}}

\numberwithin{equation}{section}

\definecolor{Maroon}{RGB}{140,10,0}

\hypersetup{
	colorlinks,
	linkcolor={Maroon},
	citecolor={Maroon},
	urlcolor={Maroon}
}

\graphicspath{{./figs/}}

\title{Convex hull of Brownian motion and Brownian bridge}

\author[S.\ \v{S}ebek]{Stjepan\ \v{S}ebek}
\address[Stjepan\ \v{S}ebek]{Department of Applied Mathematics\\
	Faculty of Electrical Engineering and Computing\\
	University of Zagreb\\
	Zagreb\\
	Croatia}
\email{stjepan.sebek@fer.hr}
\subjclass[2020]{
	60J65, % Brownian motion
	60D05, % Geometric probability and stochastic geometry
        52A22, % Random convex sets and integral geometry	
        52A10 % Convex sets in 2 dimensions
}
\keywords{convex hull, Brownian motion, Brownian bridge, area, perimeter}

\begin{document}
%\allowdisplaybreaks[4]

\begin{abstract}
In this article we study the convex hull spanned by the union of trajectories of a standard planar Brownian motion, and an independent standard planar Brownian bridge. We find exact values of the expectation of perimeter and area of such a convex hull. As an auxiliary result, that is of interest in its own right, we provide an explicit shape of the probability density function of a random variable that represents the time when combined maximum of a standard one-dimensional Brownian motion, and an independent standard one-dimensional Brownian bridge is attained. At the end, we generalize our results to the case of multiple independent standard planar Brownian motions and Brownian bridges.
\end{abstract}

\maketitle

%$\f\f\e\x\x\cdots$

%
%
% ------------------------------ INTRODUCTION ---------------------------------------
%
%
\section{Introduction}
Let $\fatW = \{\fatW(t) : t \ge 0\}$ denote the standard planar Brownian motion started at the origin, where $\fatW(t) = (W_1(t), W_2(t))$. In other words, the two coordinates $W_i(t)$, $i = 1, 2$, are independent, one-dimensional standard Brownian motions started at $0$. Furthermore, let $\fatB = \{\fatB(t) : t \in [0, 1]\}$ denote the standard planar Brownian bridge (the process that corresponds to the standard planar Brownian motion started at the origin and constrained to return to the origin at time $1$), and let us again write $\fatB(t) = (B_1(t), B_2(t))$, where the two coordinates $B_i(t)$, $i = 1, 2$, are independent, one-dimensional standard Brownian bridges started at $0$. We additionally assume that the processes $\fatW$ and $\fatB$ are independent. Denote with
\begin{align*}
    \cH^{\BM} & = \hull\{\fatW(t) : t \in [0, 1]\}, \\
    \cH^{\BB} & = \hull\{\fatB(t) : t \in [0, 1]\},
\end{align*}
and with
\begin{equation}\label{eq:def_of_H}
    \cH = \hull \left( \{\fatW(t) : t \in [0, 1]\} \cup \{\fatB(t) : t \in [0, 1]\} \right) = \hull (\cH^{\BM} \cup \cH^{\BB}),
\end{equation}
where $\hull(A)$ denotes the convex hull of $A \subseteq \bbR^2$ (the smallest convex subset of $\bbR^2$ which contains $A$). The main objects we focus on in this paper are
\begin{equation}\label{eq:def_of_L_and_A}
    \sfL = \calL(\cH), \qquad \textnormal{and} \qquad \sfA = \calA(\cH),
\end{equation}
where, $\calL(A)$ and $\calA(A)$ stand, respectively, for the perimeter and the area of the set $A \subseteq \bbR^2$. Our first main result is a closed expression for the expected value of the random variable $\sfL$. This result is formulated in the following theorem.

\begin{theorem}\label{tm:perimeter}
    Let $\sfL$ be defined as in \eqref{eq:def_of_L_and_A}. Then
    \begin{equation*}
        \bbE[\sfL] = \sqrt{2\pi} \left( 2 + \arctg \frac{1}{2} \right).
    \end{equation*}
\end{theorem}

Let us further denote $\sfL^{\BM} = \calL(\cH^{\BM})$ and $\sfL^{\BB} = \calL(\cH^{\BB})$. The problem of finding the expected value of the random variable $\sfL^{\BM}$ was proposed by Letac in 1978 \cite{Letac}, and then solved by Takacs in 1980 \cite{Takacs}. It holds that $\bbE[\sfL^{\BM}] = \sqrt{8\pi}$. When it comes to the Brownian bridge, it holds that $\bbE[\sfL^{\BB}] = \sqrt{\pi^3/2}$. This was obtained by Goldman in \cite{goldman}. 

Our next goal is to find a closed expression for the expected value of the random variable $\sfA$, defined in \eqref{eq:def_of_L_and_A}. As a step towards this goal, we developed a result which is of interest in its own right. Denote by $W = \{W(t) : t \ge 0\}$ and by $B = \{B(t) : t \in [0, 1]\}$ the standard one-dimensional Brownian motion (started at $0$), and, respectively, the standard one-dimensional Brownian bridge (started at $0$). Further, denote with $\sfT^{\BM}$ the random variable which is almost surely uniquely determined by
\begin{equation}\label{eq:def_of_TBM}
    W(\sfT^{\BM}) = \max_{0 \le s \le 1} W(s).
\end{equation}
By \cite[Theorem 2.11]{Morters-Peres} the random variable $\sfT^{\BM}$ is indeed almost surely uniquely determined by equation \eqref{eq:def_of_TBM}, and according to \cite[Theorem 5.26]{Morters-Peres} the random variable $\sfT^{\BM}$ follows the arc-sine distribution. Recall that the arc-sine distribution is a distribution on $(0, 1)$ with the probability density function $(\pi \sqrt{x(1-x)})^{-1}$, for $x \in (0, 1)$. The name of this distribution comes from the shape of its cumulative distribution function which is given by $(2 \arcsin{\sqrt{x}}) / \pi$, for $x \in (0, 1)$.

The corresponding result in the case of the standard one-dimensional Brownian bridge is also well known. Denote with $\sfT^{\BB}$ the random variable almost surely uniquely determined by
\begin{equation}\label{eq:def_of_TBB}
    B(\sfT^{\BB}) = \max_{0 \le s \le 1} B(s).
\end{equation}
It is well known that the distribution of $\sfT^{\BB}$ is uniform. This follows directly from the fact that a Brownian bridge cyclically translated an arbitrary $k \in [0, 1)$ length is still a standard Brownian bridge which has the same distribution of maximal position.

In this paper we develop an explicit formula for the probability density function of the random variable $\sfT$ which represents the time when combined maximum of a standard one-dimensional Brownian motion, and an independent standard one-dimensional Brownian bridge is attained. More precisely,
\begin{equation}\label{eq:def_of_T}
    \max\{W(\sfT), B(\sfT)\} = \max \left\{ \max_{0 \le s \le 1} W(s), \max_{0 \le s \le 1} B(s) \right\}.
\end{equation}
Notice that we do not know a priori at which of the two processes this maximum is achieved. The probability density function of the random variable $\sfT$ is given in the following theorem.

\begin{theorem}\label{tm:time_of_max}
    Let $\sfT$ be defined as in \eqref{eq:def_of_T}. Then
    \begin{equation}
        \rho(t) = \frac{4\sqrt{t}}{\pi (4t+1) \sqrt{1 - t}} + \frac{2}{\pi} \left( \frac{\sqrt{t(1-t)}}{1+t-t^2} + \arctg\sqrt{t(1-t)} \right), \quad t \in (0, 1),
    \end{equation}
    where $\rho(t)$ is the probability density function of the random variable $\sfT$.
\end{theorem}

The next result we obtain in this paper is a closed expression for the expected value of the random variable $\sfA$. This result is formulated in the following theorem.

\begin{theorem}\label{tm:area}
    Let $\sfA$ be defined as in \eqref{eq:def_of_L_and_A}. Then
    \begin{equation*}
        \bbE[\sfA] = \frac{(25 - 7\sqrt{5})\pi}{12}.
    \end{equation*}
\end{theorem}

Analogously as before, we denote $\sfA^{\BM} = \calA(\cH^{\BM})$ and $\sfA^{\BB} = \calA(\cH^{\BB})$. The expected value of $\sfA^{\BM}$ was computed by El Bachir in \cite{ElBachir_PhD_thesis}, where it is shown that $\bbE[\sfA^{\BM}] = \pi/2$. Furthermore, we have $\bbE[\sfA^{\BB}] = \pi/3$ (obtained by Majumdar et al.\ in \cite{Majumdar_EvS} and \cite{randon2009convex}). In \cite{Majumdar_EvS} the authors also give motivation for studying convex hulls of random processes, such as random walks, Brownian motion, and Brownian bridge, in terms of modelling the home range of roaming animals. This is a question that ecologists often face, in particular in designing a conservation area to preserve a given animal population \cite{murphy1992integrating}. Different methods are used to estimate this territory, based on the monitoring of the animals' positions \cite{giuggioli2006theory, worton1995convex}. One of these consists in simply the minimum convex polygon enclosing all monitored positions, that is, the convex hull. While this may seem simple minded, it remains, under certain circumstances, the best way to proceed \cite{boyle2008home}. If the recorded positions are numerous (which might result from a very fine and/or long monitoring), the number of steps of the random walker becomes large and to a good approximation the trajectory of a discrete-time planar random walk (with finite variance of the step sizes) can be replaced by a continuous-time planar Brownian motion of a certain duration. In this context, the trajectory of a Brownian bridge corresponds to an animal returning every night to its nest after spending the day foraging in the surroundings. The authors in \cite{Majumdar_EvS} then go beyond the results we mentioned so far, and study the convex hull of $n$ independent Brownian motions and, separately, the convex hull of $n$ independent Brownian bridges. However, they do not combine the two processes like we are doing here. To the best of our knowledge, the convex hull of the union of trajectories of independent standard planar Brownian motion and standard planar Brownian bridge has not been treated in the literature yet.

The rest of the paper is organized as follows. In Section \ref{sec:cauchys_formulas} we introduce the main tools we use in obtaining our results, namely Cauchy's formulas for perimeter and area. Section \ref{sec:perimeter} is devoted to developing a closed expression for the expected value of the perimeter of the convex hull $\cH$, defined in \eqref{eq:def_of_H}. In Section \ref{sec:time_of_maximum} we find the explicit formula for the probability density function of the random variable $\sfT$ defined in \eqref{eq:def_of_T} which represents the time when combined maximum of a standard one-dimensional Brownian motion, and an independent standard one-dimensional Brownian bridge is attained. In Section \ref{sec:area}, we provide a closed expression for the expected value of the area of the convex hull $\cH$. Finally, in Section \ref{sec:multiple_processes} we consider the situation when we have multiple independent copies of standard planar Brownian motion and Brownian bridge.

\section{Cauchy's formulas}\label{sec:cauchys_formulas}

The main tool that we use to find the expected value of the perimeter length and the area of the convex hull $\cH$ defined in \eqref{eq:def_of_H} are the Cauchy's integral formulas for the perimeter length and the area. The Cauchy's formula for the perimeter is a two-dimensional version of the so-called Cauchy's surface area formula (see \cite[Theorem 6.15]{Gruber_book}). An excellent reference for this formula is also \cite[Appendix A]{Majumdar_EvS} where the authors provide intuition behind this formula in two dimensions. To be able to state this formula, we need to introduce some notation. For an angle $\theta \in [0, 2\pi)$, we denote with $\bfe_{\theta}$ the unit vector in the direction of the angle $\theta$, i.e.\ $\bfe_{\theta} = (\cos \theta, \sin \theta)$. By $M(\theta)$ we denote the maximal projection of the convex hull $\cH$ onto a line passing through the origin and directed by the unit vector $\mathbf{e}_\theta$. More precisely
\begin{equation}\label{eq:def_of_M(theta)}
    M(\theta) = \max \left\{ \max_{0\le s \le 1} \langle \fatW(s), \bfe_{\theta}\rangle , \,\, \max_{0\le s \le 1} \langle\fatB(s), \bfe_{\theta}\rangle \right\},
\end{equation}
where, $\langle \cdot, \cdot \rangle$ denotes the scalar product. Cauchy's perimter formula can now be stated as
\begin{equation}\label{eq:cauchy_perimeter_formula}
    \sfL = \int_0^{2\pi} M(\theta) d\theta.
\end{equation}
Cauchy's formula for the area (we again refer the reader to \cite[Appendix A]{Majumdar_EvS}) includes the derivative of the function $M(\theta)$ with respect to $\theta$. More precisely,
\begin{equation}\label{eq:cauchy_area_formula}
    \sfA = \frac{1}{2} \int_0^{2\pi} \left( (M(\theta))^2 - (M'(\theta))^2 \right) d\theta.
\end{equation}
Directly from formulas \eqref{eq:cauchy_perimeter_formula} and \eqref{eq:cauchy_area_formula} we have
\begin{equation*}
    \bbE[\sfL] = \int_0^{2\pi} \bbE[M(\theta)] d\theta, \qquad \textnormal{and} \qquad \bbE[\sfA] = \frac{1}{2} \int_0^{2\pi} \left( \bbE[(M(\theta))^2] - \bbE[(M'(\theta))^2] \right) d\theta.
\end{equation*}
However, in the case of isotropic processes (which both standard planar Brownian motion and standard planar Brownian bridge are), these formulas become much simpler. The key consequence of the isotropic property is that the random variables $M(\theta)$ have the same law for all $\theta \in [0, 2\pi)$. Hence,
\begin{equation}\label{eq:exp_formula_for_per_and_area}
    \bbE[\sfL] = 2\pi \bbE[M(0)], \qquad \textnormal{and} \qquad \bbE[\sfA] = \pi \left(\bbE[(M(0))^2] - \bbE[(M'(0))^2]\right).
\end{equation}
It turns out that dealing with $M(0)$ is straightforward, but the random variable $M'(0)$ requires some more work. We write
\begin{equation*}
    w_{\theta}(t) = W_1(t) \cos\theta + W_2(t) \sin\theta, \quad \textnormal{and} \quad b_{\theta}(t) = B_1(t) \cos\theta + B_2(t) \sin\theta.
\end{equation*}
Taking derivatives with respect to $\theta$ we get
\begin{equation*}
        \overline{w}_{\theta}(t) = -W_1(t) \sin\theta + W_2(t) \cos\theta, \quad \textnormal{and} \quad \overline{b}_{\theta}(t) = -B_1(t) \sin\theta + B_2(t) \cos\theta.
\end{equation*}
It is clear that $w_{\theta}$ and $\overline{w}_{\theta}$ are two independent standard one-dimensional Brownian motions, and $b_{\theta}$ and $\overline{b}_{\theta}$ are two independent standard one-dimensional Brownian bridges. Hence, $M(\theta)$ (defined in \eqref{eq:def_of_M(theta)}) is simply the combined maximum of a standard one-dimensional Brownian motion $w_{\theta}$, and an independent standard one-dimensional Brownian bridge $b_{\theta}$. We now separate two scenarios, the one in which this maximum is attained by the Brownian motion $w_{\theta}$, and the one in which this maximum is attained by the Brownian bridge $b_{\theta}$. If the maximum $M(\theta)$ is attained by the Brownian motion $w_{\theta}$ at (random) time $\sfT$, we have
\begin{equation*}
    M(\theta) = w_{\theta}(\sfT) = W_1(\sfT) \cos\theta + W_2(\sfT) \sin\theta.
\end{equation*}
Therefore,
\begin{equation*}
    M'(\theta) = \overline{w}_{\theta}(\sfT) = -W_1(\sfT) \sin\theta + W_2(\sfT) \cos\theta.
\end{equation*}
In particular, $M(\theta)$ is the maximum of the first Brownian motion $w_{\theta}$, and $M'(\theta)$ corresponds to the value of the second, independent, Brownian motion $\overline{w}_{\theta}$ at time $\sfT$, when the first Brownian motion $w_{\theta}$ attains its maximum. Taking $\theta = 0$ and assuming that the combined maximum in the direction of the vector $\bfe_0$ is achieved by the Brownian motion, we have
\begin{equation}\label{eq:M_and_M'_BM_wins}
    M(0) = W_1(\sfT) = \max_{0 \le s \le 1} W_1(s), \qquad \textnormal{and} \qquad M'(0) = W_2(\sfT).
\end{equation}
Completely analogously, in the case when the combined maximum is achieved by the Brownian bridge, we have
\begin{equation}\label{eq:M_and_M'_BB_wins}
    M(0) = B_1(\sfT) = \max_{0 \le s \le 1} B_1(s), \qquad \textnormal{and} \qquad M'(0) = B_2(\sfT).
\end{equation}
For simplicity, in all the following sections we use notation $M = M(0)$, and $M' = M'(0)$.

At the end of this section we recall two well-known results that will play a crucial role in the following sections. These results give us cumulative distribution functions of random variables that represent maximum of standard one-dimensional Brownian motion, and, respectively, standard one-dimensional Brownian bridge. Using reflection principle (see \cite[Theorem 2.21]{Morters-Peres}), and, respectively, \cite[Remark 2.1]{Beghin} we get for every $x \ge 0$
\begin{equation}\label{eq:cdf_of_max_BM_and_BB}
    \bbP \left( \max_{0 \le s \le 1}W(s) \le x \right) = \erf\left( \frac{x}{\sqrt{2}} \right), \quad \textnormal{and} \quad \bbP \left( \max_{0 \le s \le 1}B(s) \le x \right) = 1 - e^{-2x^2},
\end{equation}
where
\begin{equation}\label{eq:def_of_erf}
    \erf(z) = \frac{2}{\sqrt{\pi}} \int_0^z e^{-u^2}du.
\end{equation}

\section{Expected perimeter length}\label{sec:perimeter}

In this section we prove Theorem \ref{tm:perimeter}. From equation \eqref{eq:exp_formula_for_per_and_area} it is clear that we only need to compute the expected value of the random variable $M$. Combining \eqref{eq:def_of_M(theta)}, \eqref{eq:cdf_of_max_BM_and_BB}, and the fact that $\fatW$ and $\fatB$ are independent, we have
\begin{align*}
    F_M(x)
    & = \bbP(M \le x) = \bbP \left( \max_{0 \le s \le 1}W_1(s) \le x \right) \cdot \bbP \left( \max_{0 \le s \le 1}B_1(s) \le x \right) \\
    & = \erf\left( \frac{x}{\sqrt{2}} \right) \left( 1 - e^{-2x^2} \right), \qquad x \ge 0,
\end{align*}
where we used the standard notation for the cumulative distribution function of the random variable $M$. Taking the derivative with respect to $x$ gives us the probability density function of the random variable $M$,
\begin{align}\label{eq:density_of_M}
    f_M(x)
    & = \frac{d}{dx}F_M(x) = \frac{2}{\sqrt{\pi}} e^{-\frac{x^2}{2}} \cdot \frac{1}{\sqrt{2}} \left( 1 - e^{-2x^2} \right) + \erf\left( \frac{x}{\sqrt{2}} \right) e^{-2x^2} \cdot 4x \nonumber \\
    & = \sqrt{\frac{2}{\pi}} e^{-\frac{x^2}{2}}\left( 1 - e^{-2x^2} \right) + 4xe^{-2x^2}\erf\left( \frac{x}{\sqrt{2}} \right), \qquad x \ge 0.
\end{align}
Now we have
\begin{align}\label{eq:expectation_of_M}
    \bbE[M]
    & = \int_0^{\infty} xf_M(x)dx \nonumber \\ 
    & = \sqrt{\frac{2}{\pi}} \int_0^{\infty} xe^{-\frac{x^2}{2}} \left(1 - e^{-2x^2}\right)dx + 4\int_0^{\infty} x^2 e^{-2x^2} \erf\left( \frac{x}{\sqrt{2}} \right)dx.
\end{align}
It is easy to see that
\begin{equation}\label{eq:first_int_for_exp_of_M}
    \int_0^{\infty} xe^{-\frac{x^2}{2}} \left(1 - e^{-2x^2}\right)dx = \frac{4}{5}.
\end{equation}
For the second integral in \eqref{eq:expectation_of_M} we have
\begin{align}\label{eq:second_int_for_exp_of_M}
    \int_0^{\infty}
    & x^2 e^{-2x^2} \erf\left( \frac{x}{\sqrt{2}} \right)dx = \int_0^{\infty} x^2 e^{-2x^2} \cdot \frac{2}{\sqrt{\pi}} \int_0^{x/\sqrt{2}} e^{-u^2} du dx \nonumber \\
    & = \frac{2}{\sqrt{\pi}}\int_0^{\infty} e^{-u^2} \left( \int_{u\sqrt{2}}^{\infty} x^2 e^{-2x^2} dx \right) du = \begin{bmatrix}
        u = x & dv = xe^{-2x^2} \\
        du = dx & v = -(e^{-2x^2})/4
    \end{bmatrix} \nonumber \\
    & = \frac{1}{\sqrt{2\pi}} \int_0^{\infty} ue^{-5u^2}du + \frac{1}{2\sqrt{\pi}} \int_0^{\infty}e^{-u^2} \int_{u\sqrt{2}}^{\infty} e^{-2x^2}dx du = \begin{bmatrix}
        t = x\sqrt{2} \\
        dt = \sqrt{2}dx
    \end{bmatrix} \nonumber \\
    & = \frac{1}{10\sqrt{2\pi}} + \frac{1}{2\sqrt{2\pi}} \int_0^{\infty} \int_0^{t/2} e^{-(u^2+ t^2)}du dt \nonumber \\
    & = \frac{1}{10\sqrt{2\pi}} + \frac{1}{2\sqrt{2\pi}} \int_0^{\arctg(1/2)} \int_0^{\infty}e^{-r^2} r dr d\varphi = \frac{1}{10\sqrt{2\pi}} + \frac{1}{4\sqrt{2\pi}} \arctg \frac{1}{2}.
\end{align}
Plugging \eqref{eq:first_int_for_exp_of_M} and \eqref{eq:second_int_for_exp_of_M} into \eqref{eq:expectation_of_M}, we get
\begin{equation*}
    \bbE[M] = \frac{2 + \arctg(1/2)}{\sqrt{2\pi}}.
\end{equation*}
Together with \eqref{eq:exp_formula_for_per_and_area} this gives us
\begin{equation*}
    \bbE[\sfL] = \sqrt{2\pi} \left( 2 + \arctg \frac{1}{2} \right),
\end{equation*}
which is exactly the formula from Theorem \ref{tm:perimeter}.

\section{Time of maximum}\label{sec:time_of_maximum}

As before, we denote with $\sfT$ the random variable representing the time when combined maximum of a standard one-dimensional Brownian motion, and an independent standard one-dimensional Brownian bridge is attained. In this section we find probability density function of the random variable $\sfT$, i.e.\ we prove Theorem \ref{tm:time_of_max}. Recall that the radom variable $M$ is defined as
\begin{equation*}
    M = \max \left\{\max_{0 \le s \le 1} W_1(s), \max_{0 \le s \le 1} B_1(s)\right\}.
\end{equation*}
Denote the joint probability density function of random variables $M$ and $\sfT$ with $\rho(x, t)$. In a similar fashion as before, we use notation $\rho^{\BM}(x, t)$ for the joint probability density function of the random variable representing the value of the maximum of a standard one-dimensional Brownian motion, and the time at which this maximum is achieved, and we denote with $\rho^{\BB}(x, t)$ the corresponding joint probability density function in the case of the standard one-dimensional Brownian bridge. Those joint probability density functions $\rho^{\BM}(x, t)$ and $\rho^{\BB}(x, t)$ can be computed using various techniques. As discussed in \cite[Appendix B]{Majumdar_EvS}, the simplest way is to use the Feynman-Kac path integral method, but suitably adapted with a cut-off \cite{majumdar2007brownian, majumdar2005airy}. This technique has been used \cite{majumdar2008time, randon2007distribution} to compute exactly the joint density $\rho^{\BM}(x, t)$ of a single Brownian motion, but subject to a variety of constraints, such as for a Brownian excursion, a Brownian meander etc. The results are nontrivial \cite{majumdar2008time} and have later been verified using an alternative functional renormalization group approach \cite{schehr2010extreme}. In the case of a standard one-dimensional Brownian motion, this method can be similarly used to obtain (see \cite[Appendix B]{Majumdar_EvS})
\begin{equation}\label{eq:joint_pdf_BM}
    \rho^{\BM}(x, t) = \frac{x}{\pi t^{3/2} \sqrt{1 - t}} \, e^{-\frac{x^2}{2t}}, \qquad x \ge 0, \,\, t \in (0, 1).
\end{equation}
Using the same technique one can also derive the joint probability density function of the maximum of a standard one-dimensional Brownian bridge and the time at which that maximum is achieved. It holds that (see again \cite[Appendix B]{Majumdar_EvS})
\begin{equation}\label{eq:joint_pdf_BB}
    \rho^{\BB}(x, t) = \sqrt{\frac{2}{\pi}} \, \frac{x^2}{[t(1-t)]^{3/2}} \, e^{-\frac{x^2}{2t(1-t)}}, \qquad x \ge 0, \,\, t \in (0, 1).
\end{equation}
Notice now that the combined maximum of Brownian motion and Brownian bridge will be equal to $x$ and attained at time $t$ if either Brownian motion achieves maximal value $x$ at time $t$ and maximum of Brownian bridge is less than $x$, or Brownian bridge achieves maximal value $x$ at time $t$ and maximum of Brownian motion is less than $x$. More precisely, for $x \ge 0$ and $t \in (0, 1)$, we have
\begin{equation}\label{eq:rho_0th_step}
    \rho(x, t) = \rho^{\BM}(x, t) \, \bbP\left( \max_{0\le s \le 1}B_1(s) < x \right) + \rho^{\BB}(x, t) \, \bbP\left( \max_{0\le s \le 1}W_1(s) < x \right).
\end{equation}
Combining \eqref{eq:cdf_of_max_BM_and_BB}, \eqref{eq:joint_pdf_BM} and \eqref{eq:joint_pdf_BB} we have (for $x \ge 0$ and $t \in (0, 1)$)
\begin{equation}\label{eq:rho_1st_step}
    \rho(x, t) = \frac{1}{\pi t^{3/2} \sqrt{1 - t}} \, xe^{-\frac{x^2}{2t}} \left( 1 - e^{-2x^2} \right) + \sqrt{\frac{2}{\pi}} \, \frac{1}{[t(1-t)]^{3/2}}\, x^2 e^{-\frac{x^2}{2t(1-t)}} \, \erf \left( \frac{x}{\sqrt{2}} \right).
\end{equation}
Our goal is to get the density of the random variable $\sfT$, so we have to marginalize the above joint density with respect to $x$. Using the same integration techniques as before, we get
\begin{equation}\label{eq:rho_1st_integral}
    \int_0^{\infty} xe^{-\frac{x^2}{2t}} \left( 1 - e^{-2x^2} \right)dx = \frac{4t^2}{4t+1},
\end{equation}
and furthermore
\begin{equation}\label{eq:rho_2nd_integral}
    \int_0^{\infty} x^2 e^{-\frac{x^2}{2t(1-t)}} \erf \left( \frac{x}{\sqrt{2}} \right)dx = \sqrt{\frac{2}{\pi}} t(1-t) \left[ \frac{t(1-t)}{1+t-t^2} + \sqrt{t(1-t)} \arctg \sqrt{t(1-t)} \right].
\end{equation}
Plugging \eqref{eq:rho_1st_integral} and \eqref{eq:rho_2nd_integral} into \eqref{eq:rho_1st_step} we get
\begin{equation*}
    \rho(t) = \frac{4\sqrt{t}}{\pi (4t+1) \sqrt{1 - t}} + \frac{2}{\pi} \left( \frac{\sqrt{t(1-t)}}{1+t-t^2} + \arctg\sqrt{t(1-t)} \right), \quad t \in (0, 1),
\end{equation*}
which is exactly the formula from Theorem \ref{tm:time_of_max}.

\begin{remark}\label{rem:interpretation_of_the_density}
    Notice that the part
    \begin{equation*}
        \frac{4\sqrt{t}}{\pi (4t+1) \sqrt{1 - t}}
    \end{equation*}
    corresponds to the situation when the maximum is achieved by Brownian motion. More precisely
    \begin{equation*}
        \int_0^1 \frac{4\sqrt{t}}{\pi (4t+1) \sqrt{1 - t}} dt = 1 - \frac{1}{\sqrt{5}}
    \end{equation*}
    is exactly equal to the probability that the maximum of a standard one-dimensional Brownian motion is bigger than the maximum of an independent standard one-dimensional Brownian bridge. This can also be computed in a direct way by considering the random vector $(\max_{0\le s \le 1}W_1(s), \max_{0\le s \le 1}B_1(s))$. 
\end{remark}

\section{Expected area}\label{sec:area}

In this section, we find the expected value of the random variable $\sfA$ defined in \eqref{eq:def_of_L_and_A}, and by this we prove Theorem \ref{tm:area}. We use the formula from \eqref{eq:exp_formula_for_per_and_area} which says that
\begin{equation}\label{eq:exp_formula_for_area}
    \bbE[\sfA] = \pi \left( \bbE[M^2] - \bbE[(M')^2] \right).
\end{equation}
Clearly, we need to compute $\bbE[M^2]$ and $\bbE[(M')^2]$. We first compute $\bbE[M^2]$. Using \eqref{eq:density_of_M}, we have
\begin{equation}\label{eq:expectation_of_M^2_1st_step}
    \bbE[M^2] = \sqrt{\frac{2}{\pi}} \int_0^{\infty} x^2 e^{-\frac{x^2}{2}} \left( 1 - e^{-2x^2} \right)dx + 4\int_0^{\infty} x^3 e^{-2x^2} \erf\left( \frac{x}{\sqrt{2}} \right)dx.
\end{equation}
It is easy to see that
\begin{equation}\label{eq:first_int_for_exp_of_M^2}
    \int_0^{\infty} x^2 e^{-\frac{x^2}{2}} \left( 1 - e^{-2x^2} \right)dx = \frac{\sqrt{2\pi}(5\sqrt{5} - 1)}{10\sqrt{5}},
\end{equation}
and
\begin{equation}\label{eq:second_int_for_exp_of_M^2}
    \int_0^{\infty} x^3 e^{-2x^2} \erf\left( \frac{x}{\sqrt{2}} \right)dx = \frac{7}{40\sqrt{5}}.
\end{equation}
Plugging \eqref{eq:first_int_for_exp_of_M^2} and \eqref{eq:second_int_for_exp_of_M^2} into \eqref{eq:expectation_of_M^2_1st_step} we get
\begin{equation}\label{eq:expectation_of_M^2}
    \bbE[M^2] = 1 + \frac{1}{2\sqrt{5}}.
\end{equation}
The last part is to compute $\bbE[(M')^2]$. From \eqref{eq:M_and_M'_BM_wins} and \eqref{eq:M_and_M'_BB_wins} we know that $M'$ corresponds either to $W_2(\sfT)$ (where $\sfT$ is the time when the maximum is achieved, and the maximum is achieved by Brownian motion) or to $B_2(\sfT)$ (when the maximum is achieved by Brownian bridge). For simplicity, denote by
\begin{equation*}
    A = \left\{ \max_{0 \le s \le 1}W_1(s) > \max_{0 \le s \le 1}B_1(s) \right\},
\end{equation*}
and by $A^c$ the complement of the set $A$. Using well-known facts that $\bbE[(W_2(t))^2] = t$ and $\bbE[(B_2(t))^2] = t(1-t)$ (see also Remark \ref{rem:interpretation_of_the_density}), we have
\begin{equation}\label{eq:calculation_of_exp_of_M'2}
    \begin{aligned}
        \bbE[(M')^2]
        & = \bbE[(M')^2 \bbjedan_A] + \bbE[(M')^2 \bbjedan_{A^c}] \\
        & = \bbE[(W_2(\sfT))^2 \bbjedan_A] + \bbE[(B_2(\sfT))^2 \bbjedan_{A^c}] \\
        & = \int_0^1 \int_0^{\infty} \bbE[(W_2(t))^2] \rho^{\BM}(x, t) \bbP \left( \max_{0 \le s \le 1}B_1(s) < x\right) dx dt \\
        & \qquad + \int_0^1 \int_0^{\infty} \bbE[(B_2(t))^2] \rho^{\BB}(x, t) \bbP \left( \max_{0 \le s \le 1}W_1(s) < x\right) dx dt \\
        & = \int_0^1 t \cdot \frac{4\sqrt{t}}{\pi (4t+1) \sqrt{1-t}} dt \\
        & \qquad + \int_0^1 t(1-t) \cdot \frac{2}{\pi} \left[ \frac{\sqrt{t(1-t)}}{1+t-t^2} + \arctg \sqrt{t(1-t)} \right]dt \\
        & = \frac{4}{\pi} \int_0^1 \frac{t\sqrt{t}}{(4t+1)\sqrt{1-t}}dt + \frac{2}{\pi}\int_0^1 \frac{[t(1-t)]^{3/2}}{1+t-t^2}dt \\
        & \qquad + \frac{2}{\pi} \int_0^1 t(1-t) \arctg \sqrt{t(1-t)}dt \\
        & = \frac{4}{\pi} \cdot \frac{(5+\sqrt{5})\pi}{80} + \frac{2}{\pi}\cdot \frac{(16 - 7\sqrt{5})\pi}{8\sqrt{5}} + \frac{2}{\pi} \cdot \frac{(5\sqrt{5} - 10)\pi}{24\sqrt{5}} \\
        & = \frac{41 - 13\sqrt{5}}{12\sqrt{5}}.
    \end{aligned}
\end{equation}
Plugging this and \eqref{eq:expectation_of_M^2} into formula \eqref{eq:exp_formula_for_area} gives us
\begin{equation*}
    \bbE[\sfA] = \pi \left( 1 + \frac{1}{2\sqrt{5}} - \frac{41 - 13\sqrt{5}}{12\sqrt{5}} \right) = \frac{(25 - 7\sqrt{5})\pi}{12},
\end{equation*}
which is exactly the formula from Theorem \ref{tm:area}.

\section{Multiple processes}\label{sec:multiple_processes}

In this section we generalize our main results to the case of multiple independent copies of standard planar Brownian motion and Brownian bridge. Let $m, n \in \bbN$ be arbitrary positive integers. Let $\fatW^1, \fatW^2, \ldots, \fatW^m$ (where $\fatW^i = (W_1^i, W_2^i)$) be $m$ independent copies of the process $\fatW$, and $\fatB^1, \fatB^2, \ldots, \fatB^n$ (where $\fatB^j = (B_1^j, B_2^j)$) be $n$ independent copies of the process $\fatB$. We assume further that $\fatW^i$ and $\fatB^j$ are independent for all $i \in \{1, 2, \ldots, m\}$ and $j \in \{ 1, 2, \ldots, n\}$. Due to isotropic property of both considered processes, we can again simplify formulas for computing the expected perimeter and area. Denote the perimeter and the area of the convex hull spanned by $m$ independent standard planar Brownian motions run up to time $1$, and $n$ independent standard planar Brownian bridges with $\sfL_{m, n}$, and $\sfA_{m, n}$ respectively. Analogously as in
\eqref{eq:exp_formula_for_per_and_area}, we have
\begin{equation}\label{eq:exp_formula_per_and_ar_multiple_proc}
    \bbE[\sfL_{m, n}] = 2\pi \bbE[M_{m, n}], \qquad \textnormal{and} \qquad \bbE[\sfA_{m, n}] = \pi\left( \bbE[M_{m, n}^2] - \bbE[(M_{m, n}')^2] \right),
\end{equation}
where
\begin{equation}\label{eq:def_of_Mmn}
    M_{m, n} = \max\left\{ \max_{1 \le i \le m} \max_{0 \le s \le 1} W_1^i(s), \max_{1 \le j \le n} \max_{0 \le s \le 1} B_1^j(s) \right\},
\end{equation}
and
\begin{equation}\label{eq:def_of_Mmn'}
    M_{m, n}' = \begin{cases}
        W_i^2(\sfT_{m, n}), & \textnormal{if max is attained by } i \textnormal{-th Brownian motion at time } \sfT_{m, n}, \\
        B_j^2(\sfT_{m, n}), & \textnormal{if max is attained by } j \textnormal{-th Brownian bridge at time } \sfT_{m, n}.
    \end{cases}
\end{equation}
Let us first find the probability density function of the random variable $M_{m, n}$. From \eqref{eq:def_of_Mmn} and \eqref{eq:cdf_of_max_BM_and_BB} we have
\begin{equation*}
    F_{M_{m, n}}(x) = \bbP(M_{m, n} \le x) = \erf\left( \frac{x}{\sqrt{2}} \right)^m \left( 1 - e^{-2x^2} \right)^n, \qquad x \ge 0.
\end{equation*}
From this, we directly have
\begin{equation}\label{eq:pdf_of_Mmn}
    \begin{aligned}
        f_{M_{m, n}}(x) = \frac{d}{dx} F_{M_{m, n}}(x)
        & = \sqrt{\frac{2}{\pi}} m e^{-\frac{x^2}{2}} \erf\left( \frac{x}{\sqrt{2}} \right)^{m - 1} \left( 1 - e^{-2x^2} \right)^n \\
        & \qquad + 4nx e^{-2x^2} \erf\left( \frac{x}{\sqrt{2}} \right)^m \left( 1 - e^{-2x^2} \right)^{n-1}, \qquad x \ge 0.
    \end{aligned}
\end{equation}
Combining this with \eqref{eq:exp_formula_per_and_ar_multiple_proc} gives us
\begin{equation}\label{eq:formula_for_Lmn}
    \begin{aligned}
        \bbE[\sfL_{m, n}]
        & = 2m \sqrt{2\pi} \int_0^{\infty} x e^{-\frac{x^2}{2}} \erf\left( \frac{x}{\sqrt{2}} \right)^{m - 1} \left( 1 - e^{-2x^2} \right)^n dx \\
        & \qquad + 8n\pi \int_0^{\infty} x^2 e^{-2x^2} \erf\left( \frac{x}{\sqrt{2}} \right)^m \left( 1 - e^{-2x^2} \right)^{n-1} dx.
    \end{aligned}
\end{equation}
Plugging $m = n = 1$ we reconstruct the result from Theorem \ref{tm:perimeter}. The next few values are
\begin{align*}
    \bbE[\sfL_{1, 2}] & = \frac{64\sqrt{2\pi}}{45} + \frac{\sqrt{\pi}}{45} \left( 26\sqrt{2} + 90\sqrt{2} \arcctg(2) - 45 \arcctg(2\sqrt{2}) \right) \approx 6.7353,\\
    \bbE[\sfL_{2, 1}] & \approx 7.5945,\\
    \bbE[\sfL_{2, 2}] & \approx 7.9019.
\end{align*}
While performing the calculation of these expected values, one can see that some of the integrals appearing can be computed explicitly, but some we only managed to compute numerically. As soon as there was at least one integral appearing in the calculations that we only managed to evaluate numerically, we presented the numerical value for the final result.

We use notation $\sfT_{m, n}$ for the random variable representing the time when combined maximum of $m$ independent standard one-dimensional Brownian motions ($W_1^1, W_1^2, \ldots, W_1^m$), and $n$ independent standard one-dimensional Brownian bridges ($B_1^1, B_1^2, \ldots, B_1^n$) is attained. Analogously as in \eqref{eq:rho_0th_step}, using \eqref{eq:cdf_of_max_BM_and_BB}, \eqref{eq:joint_pdf_BM} and \eqref{eq:joint_pdf_BB} we have
\begin{align*}
    \rho_{m, n}(x, t)
    & = m \rho^{\BM}(x, t) \bbP \left( \max_{0 \le s \le 1} W_1^1(s) < x \right)^{m - 1} \bbP \left( \max_{0 \le s \le 1} B_1^1(s) < x \right)^n \\
    & \qquad + n \rho^{\BB}(x, t) \bbP \left( \max_{0 \le s \le 1} W_1^1(s) < x \right)^m \bbP \left( \max_{0 \le s \le 1} B_1^1(s) < x \right)^{n-1} \\
    & = \frac{m}{\pi t^{3/2}\sqrt{1-t}} \, xe^{-\frac{x^2}{2t}}\erf\left( \frac{x}{\sqrt{2}} \right)^{m-1}  \left( 1-e^{-2x^2} \right)^n \\
    & \qquad + \sqrt{\frac{2}{\pi}} \frac{n}{[t(1-t)]^{3/2}} \, x^2 e^{-\frac{x^2}{2t(1-t)}} \erf\left( \frac{x}{\sqrt{2}} \right)^m  \left( 1-e^{-2x^2} \right)^{n-1},
\end{align*}
where $\rho_{m, n}(x,t)$ is the joint probability density function of random variables $M_{m, n}$ and $\sfT_{m, n}$. Denote with
\begin{align*}
    \rho_{m, n}^1(t) & = \frac{m}{\pi t^{3/2}\sqrt{1-t}} \int_0^{\infty} xe^{-\frac{x^2}{2t}}\erf\left( \frac{x}{\sqrt{2}} \right)^{m-1}  \left( 1-e^{-2x^2} \right)^n dx, \\
    \rho_{m, n}^2(t) & = \sqrt{\frac{2}{\pi}} \frac{n}{[t(1-t)]^{3/2}} \int_0^{\infty} x^2 e^{-\frac{x^2}{2t(1-t)}} \erf\left( \frac{x}{\sqrt{2}} \right)^m  \left( 1-e^{-2x^2} \right)^{n-1} dx.
\end{align*}
Denoting the density of the random variable $\sfT_{m, n}$ with $\rho_{m, n}(t)$, we have $\rho_{m, n}(t) = \rho_{m, n}^1(t) + \rho_{m, n}^2(t)$. Analogously as in \eqref{eq:calculation_of_exp_of_M'2}, we have
\begin{equation}\label{eq:exp_of_Mmn'_with_rho}
    \bbE[(M_{m, n}')^2] = \int_0^1 t \rho_{m, n}^1(t) dt + \int_0^1 t(1-t) \rho_{m, n}^2(t)dt.
\end{equation}
Combining \eqref{eq:exp_formula_per_and_ar_multiple_proc}, \eqref{eq:pdf_of_Mmn} and \eqref{eq:exp_of_Mmn'_with_rho} we get
\begin{align*}
    \bbE[\sfA_{m, n}]
        & = m \sqrt{2\pi} \int_0^{\infty} x^2 e^{-\frac{x^2}{2}} \erf\left( \frac{x}{\sqrt{2}} \right)^{m - 1} \left( 1 - e^{-2x^2} \right)^n dx \\
        & \qquad + 4n\pi \int_0^{\infty} x^3 e^{-2x^2} \erf\left( \frac{x}{\sqrt{2}} \right)^m \left( 1 - e^{-2x^2} \right)^{n-1} dx \\
        & \qquad \qquad - \pi \left( \int_0^1 t \rho_{m, n}^1(t) dt + \int_0^1 t(1-t) \rho_{m, n}^2(t)dt \right).
\end{align*}
Plugging $m = n = 1$ we reconstruct the result from Theorem \ref{tm:area}. The next few values are
\begin{align*}
    \bbE[\sfA_{1, 2}] & \approx 2.9705 ,\\
    \bbE[\sfA_{2, 1}] & \approx 3.6966 ,\\
    \bbE[\sfA_{2, 2}] & \approx 4.0651.
\end{align*}
As in the case of the perimeter, some of the integrals appearing can be computed explicitly, but in every expression there is at least one integral that we only managed to compute numerically.

\begin{remark}
    In \cite{Majumdar_EvS}, the authors provide heuristic argument that the convex hull of $m$ independent standard planar Brownian motions run up to time $1$ (for large $m$) resembles the circle centered at the origin, with radius $\sqrt{2\ln m}$. More precisely, the Hausdorff distance between the convex hull of $m$ independent standard planar Brownian motions run up to time $1$, and the circle centered at the origin, with radius $\sqrt{2\ln m}$, converges to $0$, as $m$ tends to infinity. Analogously, they argue that the Hausdorff distance between the convex hull of $n$ independent standard planar Brownian bridges, and the circle centered at the origin, with radius $\sqrt{(\ln n) / 2}$, converges to $0$, as $n$ tends to infinity. These observations were made formal in \cite{davydov2011convex}. From this, it is clear that the convex hull spanned by $m$ independent standard planar Brownian motions and $n$ independent standard planar Brownian bridges also approaches a circle. If $n < m^4$, the radius of that circle will be dictated with the number of Brownian motions, and if $n > m^4$, the convex hull generated by Brownian bridges wins.
\end{remark}

\section*{Acknowledgments}
\noindent
Financial support of the Croatian Science Foundation (project IP-2022-10-2277) is gratefully acknowledged.

\bibliographystyle{bababbrv-fl}
\bibliography{literature}

\end{document}